\documentclass[10pt]{article}
\usepackage[a4paper, top=2cm, bottom=2cm, left=2cm, right=2.5cm]{geometry}
\usepackage{mathtools}
\usepackage{amsmath}
\usepackage{caption}
\title{\bfseries A Vector-Algebraic Reconstruction of the Medieval Arab Formula\\
for the Shadow of a Gnomon}

\author{Fabrizio Patuzzo\\[4pt]
\small Independent Researcher\\
\small \texttt{fabrizio.patuzzo@gmail.com}}

\date{Oct 19, 2025} 

\begin{document}

\maketitle
\thispagestyle{plain}

\begin{abstract}
We present a modern reconstruction of the classical formula—first derived by medieval Arab astronomers—that describes the trajectory of the tip of a gnomon’s shadow during the day as a function of latitude, solar declination, and gnomon height. Unlike the traditional derivations based on spherical trigonometry, our approach uses only elementary vector algebra and rotation matrices. The same vector framework naturally extends to other classical relations used in sundial construction.
\end{abstract}

\section{Introduction}
The mathematical analysis of sundials, or gnomonics, dates back to antiquity and reached a high level of sophistication among medieval Arab astronomers. They developed explicit formulas describing the position of the shadow cast by a vertical gnomon as the Sun moves across the sky.
The ‘Arab formula’ appears in various treatises from the 9th to the 13th centuries, such as those of al-Battānī and Ibn al-Shāṭir, and was traditionally obtained through spherical trigonometry, often with intricate geometric reasoning.\\[6pt]
In this paper we show that the same results can be derived directly from elementary vector algebra, using successive rotations of the vectors of the Sun’s direction, the gnomon and the dial plane. This approach highlights the underlying simplicity of the geometry and provides an intuitive modern reformulation of the classical results.

\section{The Arab formula}
The mathematical expression known as the Arab formula describes the trajectory of the tip of the shadow cast by a vertical gnomon during the course of a day. For a gnomon of height  $h$, placed at latitude $\lambda$, and for a solar declination $\delta$, the locus of the shadow tip satisfies \\[6pt]
\begin{equation}
y = \frac{-h \cdot sin\lambda \cdot cos\lambda + sin\delta \sqrt{(cos^2 \lambda - sin^2 \delta)x^2 +h^2 cos^2 \delta}}{sin^2 \delta - cos^2\lambda}
\label{eq:arab_formula}
\end{equation}
At first glance this formula may appear forbiddingly complex—“it looks like Arabic,” one might say—in fact it was precisely Arab astronomers of the Middle Ages who first derived it. In what follows we reconstruct their result step by step, showing that it can be obtained using nothing more than elementary vector algebra.

\subsection{A Simple Case: the North Pole at the Equinox}
We begin with the simplest possible configuration: an observer located at the North Pole ($\lambda=90^{\circ}$) at the equinox ($\delta=0$). At this moment the Sun lies on the celestial equator and its rays reach the observer horizontally. A vertical gnomon, represented by the vector $\mathbf{g}=(0,0,h)$, is therefore perpendicular to the Sun’s direction $\mathbf{s}=(0,1,0)$ and the angle between them is simply $\alpha=\arccos{(\mathbf{g} \cdot \mathbf{s})} = 90^{\circ}$. The computation is trivial, but it illustrates the key idea: the geometry of the shadow can be obtained once the two direction vectors—the Sun’s rays and the gnomon—are known.

\subsection{The Sun’s Direction Vector}
To describe the Sun’s direction at any date and time, we begin from its equinoctial position and apply successive rotations.

\begin{enumerate}
\item Rotation by the declination
$\delta$ around the $x$-axis:
\begin{equation}
\mathbf{s}_1=\begin{bmatrix}
1 & 0 & 0\\
0 & \cos\delta & \sin\delta\\
0 & -\sin\delta & \cos\delta
\end{bmatrix} \cdot \begin{bmatrix}
0 \\ 1 \\ 0
\end{bmatrix}=\begin{bmatrix}
 0 \\
\cos\delta \\ 
-\sin\delta
\end{bmatrix}
\end{equation}

\item Rotation by the hour angle
H around the z-axis:
\begin{equation}
\mathbf{s}_2=\begin{bmatrix}
\cos H & \sin H & 0\\
-\sin H & \cos H & 0\\
0 & 0 & 1
\end{bmatrix} \cdot \mathbf{s}_1=\begin{bmatrix}
\cos\delta \cdot \sin H \\
\cos\delta \cdot \cos H \\
-\sin\delta
\end{bmatrix}
\end{equation}

\end{enumerate}
The vector $\mathbf{s}_2$ gives the direction of the solar rays for a given declination and hour angle.

\subsection{The Arab Formula when $\lambda=90^{\circ}$}
One can find the formula describing the trajectory of the tip of the gnomon's shadow in a given day at the North Pole by solving the equation
\begin{equation}
q\cdot \begin{bmatrix}
\cos\delta \cdot \sin H \\
\cos\delta \cdot \cos H \\
-\sin\delta
\end{bmatrix}+\begin{bmatrix}
0 \\ 0 \\ h
\end{bmatrix}=\begin{bmatrix}
x \\ y \\ 0
\end{bmatrix}
\end{equation}
which yields 
\begin{equation}
x^2+y^2=(h \cdot \cot\delta)^2
\end{equation}
The trajectory is a circle of radius $h \cdot \cot\delta$, consistent with the continuous rotation of the Sun around the horizon at the Pole, a result that one can also obtain by setting $\lambda=90^{\circ}$ in Eq.~\eqref{eq:arab_formula}.

\subsection{The Gnomon and Local Horizontal Plane}
At latitude $\lambda$, the vertical gnomon is no longer perpendicular to the equatorial plane. We represent this by rotating the gnomon vector by $(90-\lambda)^{\circ}$ around the $x$-axis, obtaining
\begin{equation}
\mathbf{g}_1=
\begin{bmatrix}
1 & 0 & 0\\
0 & \sin\lambda & -\cos\lambda\\
0 & \cos\lambda & \sin\lambda
\end{bmatrix}\cdot 
\begin{bmatrix}
0\\
0\\
h
\end{bmatrix}
=
\begin{bmatrix}
0\\
-h\cdot \cos\lambda\\
h \cdot \sin\lambda
\end{bmatrix}
\end{equation}
Last, we represent the coordinates of points on the local horizontal plane at latitude $\lambda$ with the vector
\begin{equation}
\mathbf{z}_1=\begin{bmatrix}
1 & 0 & 0\\
0 & \sin\lambda & -\cos\lambda\\
0 & \cos\lambda & \sin\lambda
\end{bmatrix} \cdot \begin{bmatrix}
x\\
y\\
0
\end{bmatrix}=\begin{bmatrix}
x\\
y\cdot sin\lambda\\
y\cdot cos\lambda
\end{bmatrix}
\end{equation}



\subsection{The Arab Formula via Vector Algebra}
The following equation describes the intersection between the line representing the sun's ray passing through the tip of the gnomon and the horizontal plane at latitude $\lambda$ \\[6pt]
\begin{equation}
q\cdot \begin{bmatrix}
\cos\delta \cdot \sin H \\
\cos\delta \cdot \cos H \\
-\sin\delta
\end{bmatrix} + \begin{bmatrix}
0\\
-h\cdot cos\lambda\\
h \cdot sin\lambda
\end{bmatrix} = 
\begin{bmatrix}
x\\
y\cdot sin\lambda\\
y\cdot cos\lambda
\end{bmatrix}
\end{equation}

By solving for $y$, we recover Eq.~\eqref{eq:arab_formula} 
\begin{equation}
y = \frac{-h \cdot sin\lambda \cdot cos\lambda + sin\delta \sqrt{(cos^2 \lambda - sin^2 \delta)x^2 +h^2 cos^2 \delta}}{sin^2 \delta - cos^2\lambda}
\end{equation}
the historical Arab formula.

\section{Interpretation and Conic Classification}

Eq.~\eqref{eq:arab_formula} can be interpreted geometrically as the intersection
between two simple surfaces in three--dimensional space:
the cone of solar rays and the horizontal plane of the ground.
The conic nature of the shadow’s trajectory therefore follows immediately
from elementary geometry.

\subsection{From the vector equation to the conic form}

The tip of the gnomon defines the vertex of a cone whose
generating lines are parallel to the Sun’s direction vector~$\mathbf{s}_2$. Its intersection with the ground plane~$z=0$
gives the curve described by the shadow during the day. This expression may be rewritten in the standard quadratic form
\begin{equation}
A x^{2} + B x y + C y^{2} + D x + E y + F = 0,
\label{eq:general_conic}
\end{equation}
whose coefficients depend on the observer’s latitude~$\lambda$,
the solar declination~$\delta$,
and the gnomon height~$h$.

\subsection{Classification of the shadow path}

The nature of the curve can be determined from the discriminant
$\Delta = B^{2} - 4AC$ of Eq.~\eqref{eq:general_conic}.
A straightforward computation shows that
\begin{equation}
\Delta \propto  \cos^{2}\lambda-\sin^{2}\delta.
\end{equation}
Hence:
\[
\begin{cases}
\cos^{2}\lambda - \sin^{2}\delta>0  &\Rightarrow \text{hyperbola,}\\[4pt]
\cos^{2}\lambda - \sin^{2}\delta<0  &\Rightarrow \text{ellipse or circle},\\[4pt]
\cos^{2}\lambda -\sin^{2}\delta = 0 &\Rightarrow \text{parabola}.
\end{cases}
\]
The daily path of the shadow forms a hyperbola at low and middle latitudes, an ellipse or circle at polar latitudes, and a parabola in the limiting case when $\lambda=90^{\circ}-\delta$, for example at the Arctic circle during the summer solstice.

\section{Relation to Other Sundial Formulas}

The vector formulation developed above is not limited to the derivation of the shadow--tip trajectory. It also allows one to obtain many of the classical formulas used in sundial construction with little additional effort.

\subsection{Sunrise and sunset times}

The times of sunrise and sunset happen when the sun's direction vector belongs to the local horizontal plane at latitude $\lambda$. Hence 
\begin{equation}
\begin{bmatrix}
cos\delta \cdot sinH_0 \\
cos\delta \cdot cosH_0 \\
-sin\delta
\end{bmatrix}= 
\begin{bmatrix}
x \\
y \cdot sin\lambda \\
y \cdot cos\lambda
\end{bmatrix}
\end{equation}
which yields
\begin{equation}
\cos H_0 = -\tan\lambda \tan\delta,
\end{equation}
The sunrise and sunset times are
\begin{equation}
sunrise = 12h - \frac{\arccos (- \tan\lambda \tan\delta)}{15}
\end{equation}
\begin{equation}
sunset = 12h + \frac{\arccos (- \tan\lambda \tan\delta)}{15}
\end{equation}
where the angles are expressed in degrees.

\subsection{Solar altitude}

The altitude $a$ of the Sun above the horizon can be obtained from the dot product between the vector representing the sun direction and the vector representing a gnomon of unit height
\begin{equation}
\sin a = -\begin{bmatrix}
\cos\delta \cdot \sin H \\
\cos\delta \cdot \cos H \\
-\sin\delta
\end{bmatrix} \cdot \begin{bmatrix}
0\\
-\cos\lambda\\
\sin\lambda
\end{bmatrix} 
\end{equation}
which yields
\begin{equation}
\sin a = \sin\lambda \sin\delta + \cos\lambda \cos\delta \cos H
\end{equation}

\subsection{Hour angle}

To find the hour angle $H'$ of a polar gnomon at latitude $\lambda$ we first solve the equation
\begin{equation}
q\cdot \begin{bmatrix}
\cos\delta \cdot \sin H \\
\cos\delta \cdot \cos H \\
-\sin\delta
\end{bmatrix} + \begin{bmatrix}
0\\
0\\
h 
\end{bmatrix} = 
\begin{bmatrix}
x\\
y\cdot sin\lambda\\
y\cdot cos\lambda
\end{bmatrix}
\end{equation}

and then divide $x$ by $y$, which yields
\begin{equation}
tanH'=tanH \cdot sin\lambda
\end{equation}

\section{Conclusion}

The derivation presented in this paper shows that the classical ``Arab formula''  describing the path of a gnomon's shadow
can be reconstructed entirely within the framework of elementary vector algebra.\\[6pt]
This approach also clarifies the underlying structure of other classical sundial relations, such as those governing
sunrise and sunset times, solar altitude, and hour lines. \\[6pt]
The method offers pedagogical value: it demonstrates how a few rotations and dot products can reproduce the results that once required pages of spherical geometry.
The same ideas can be extended to the construction of general sundial surfaces. The ancient art of sundials continues to illuminate the teaching of geometry.

\end{document}